\documentclass[12pt]{article}
\usepackage[intlimits]{amsmath}
\usepackage{amsfonts,amssymb,amscd,amsthm}
\usepackage{psi-d-o}

\newcommand{\Psiw}{\Psi_\wedge^0}

\DeclareMathOperator{\Mat}{Mat}
\DeclareMathOperator{\Con}{Con}

\newtheorem{theorem}{Theorem}[section]
\newtheorem{corollary}[theorem]{Corollary}
\newtheorem{lemma}[theorem]{Lemma}
\newtheorem{proposition}[theorem]{Proposition}

\theoremstyle{definition}

\newtheorem{remark}[theorem]{Remark}
\newtheorem{example}[theorem]{Example}
\newtheorem{definition}[theorem]{Definition}

\newlength{\rig}
\newlength{\rigg}
\newlength{\hei}

\numberwithin{equation}{section}

\title{\bf On the Homotopy Classification of Elliptic Operators on Manifolds
with Edges}

\author{{\bf V.E.~Nazaikinskii, A.~Yu.~Savin, B.-W.~Schulze, B.~Yu.~Sternin}}

\date{}

\begin{document}
\maketitle

\begin{abstract}
We obtain a stable homotopy classification of elliptic operators on manifolds with
edges.
\end{abstract}

\tableofcontents

\section*{Introduction} \addcontentsline{toc}{section}{Introduction}

The paper deals with the classification of elliptic operators on manifolds with edges, i.e., a
description of the set of elliptic operators up to stable homotopy.

For the first time, such a classification for the case of smooth manifolds was given by Atiyah
and Singer \cite{AtSi1} in terms of the topological $K$-functor. Later, an approach to the
homotopy classification on manifolds with singularities in terms of the \emph{analytic
$K$-homology} of the manifold was suggested in~\cite{Sav8}. Namely, an elliptic operator on a
manifold $\mathcal{M}$ with singularities is represented as an abstract (elliptic) operator in
the sense of Atiyah \cite{Ati4} and hence defines a cycle in the analytic $K$-homology of the
manifold $\cM$ viewed as a compact topological space.  In other words, there is a well-defined
homomorphism
\begin{equation}\label{klass}
\Ell_0(\mathcal{M})\longrightarrow  K_0(\mathcal{M}),
\end{equation}
where $\Ell_0(\mathcal{M})$ is the group of stable homotopy classes of elliptic operators. For
the case of isolated singularities, the fact that \eqref{klass} is an isomorphism was proved
in~\cite{Sav8}. (See also~\cite{NScS5}, where the case of one singular point was considered.)

The classification \eqref{klass} has many corollaries and applications, including a formula for
the obstruction, similar to the Atiyah--Bott obstruction~\cite{AtBo2}, to the existence of
Fredholm problems for elliptic equations, the fact that the group $\Ell_0(\mathcal{M})$ is
equal modulo torsion to the homology $H_{ev}(\mathcal{M})$, a generalization of Poincar\'e
duality in $K$-theory to manifolds with singularities, etc.

In the present paper, we prove that the mapping \eqref{klass} is an isomorphism for elliptic
operators on manifolds with edges in the sense of \cite{NSScS98}.

The idea underlying the proof is simple. A manifold $\mathcal{M}$ with edges is a stratified
manifold with two strata, the singularity stratum $X$ and the open stratum
$\mathcal{M}\setminus X$; both strata are smooth. The isomorphism \eqref{klass} on smooth
manifolds is known \cite{BaDo1,Kas3}, and hence it is natural to extend the assertion about the
isomorphism to the union of these two strata using the sequences
\begin{equation}
\label{fantomth}
\begin{array}{ccccccccc}
     \Ell_1(\mathcal{M}\setminus X) &\to &    \Ell_0(X) & \to&
\Ell_0(\mathcal{M})& \to &\Ell_0(\mathcal{M}\setminus X) &\to
&\Ell_1(X) \\
\downarrow & & \downarrow & & \downarrow & & \downarrow & & \downarrow  \\
K_1(\mathcal{M}\setminus X) &\to &    K_0(X) & \to& K_0(\mathcal{M})& \to
&K_0(\mathcal{M}\setminus X) &\to &K_1(X)
\end{array}
\end{equation}
of the pair $X\subset \mathcal{M}$ in $K$-homology and $\Ell$-theory. Once this commutative
diagram is defined, the desired isomorphism follows from the five lemma.

Actually, we construct $\Ell$-theory ($\Ell$-functor) for the manifold $\mathcal{M}$ with edges
and then establish the isomorphism \eqref{klass} by mimicking the proof of the uniqueness
theorem for extraordinary cohomology theory. We however point out the following important
facts.
\begin{enumerate}
\item Exact sequences in elliptic theory are not known in general. To construct these exact
sequences, we represent the group $\Ell_0(\mathcal{M})$ as the $K$-group of some $C^*$-algebra
and then use the exact sequence of $K$-theory of algebras.

\item The main difficulty is the boundary map in the upper row (i.e., the boundary map in
$K$-theory of $C^*$-algebras). We use \textit{semiclassical quantization} (e.g.,
see~\cite{NScS12}), which permits us to replace the algebra of edge symbols by a simpler
algebra of families of parameter-dependent operators in the computation of the boundary map.
For the latter algebra, the boundary map is given by the index theorem in~\cite{Nis1}.
\end{enumerate}

Let us briefly outline the contents of the paper. We recall the main notions of elliptic theory
on manifolds with edges in Sec.~1. Section~2 describes a construction that assigns a cycle in
$K$-homology to each elliptic operator. Then (Sec.~3) we state a homotopy classification
theorem and present its proof except for two especially lengthy computations, which are given
separately in Secs.~4 and~5. The last section contains some additional remarks (the
classification of edge morphisms and the topological obstruction to the existence of elliptic
edge problems). The desired properties of semiclassical quantization are established in the
Appendix.

The research was supported by the Deutsche Forschungsgemeinschaft  and by RFBR grants
Nos.~02-01-00118, 02-01-00928, and~03-02-16336.

\section{Operators on Manifolds with Edges}

First, we recall some facts of the theory of elliptic operators on manifolds with edges. We
systematically use the results of the paper \cite{NSScS98}.\footnote{The paper \cite{NSScS98}
contains a version of the theory of elliptic operators on manifolds with edges
\cite{Schu1,EgSc1} especially suited for the study of topological aspects of the theory. In
particular, one deals only with smooth (or continuous) symbols and does not impose any
analyticity requirements.}

\paragraph{1. Manifolds with edges.} Let $M$ be a smooth compact manifold with boundary
$\partial M$ equipped with the structure of a smooth locally trivial bundle $\pi:\partial M \to
X$ with base $X$ and fiber $\Omega$. A \textit{manifold $\mathcal{M}$ with edge} $X\subset
\mathcal{M}$ is the space obtained from $M$ by identifying the points lying in the same fiber:
$$
\mathcal{M}=M\bigr/ \sim,\qquad x\sim y \Longleftrightarrow x=y, \text{ or } (x,y\in\partial M
\text{ and } \pi(x)=\pi(y)).
$$
The complement $\ciM=\mathcal{M}\setminus X$ is an open smooth manifold, and an arbitrary point
of the edge $X$ has a neighborhood homeomorphic to the \emph{model wedge}
\begin{equation}\label{modell}
\WW=\RR^n\times K_\Omega,
\end{equation}
where $n$ is the dimension of $X$ and $K_\Omega=\Omega\times
\overline{\mathbb{R}}_+\bigr/\Omega\times \{0\}$ is the cone with base $\Omega$. In this local
model, the points of the edge form the subset  $\mathbb{R}^n\times
\{0\}\subset\mathbb{R}^n\times K_\Omega$.

We often use local coordinates $(x,\omega,r)$.  For the model wedge \eqref{modell}, these
coordinates have the following form: $x\in \mathbb{R}^n$, $r\in {\mathbb{R}}_+$, and $\omega$
is a coordinate on $\Omega$.

\paragraph{2. Differential operators and function spaces.}
Consider a differential expression with smooth coefficients on $\ciM$ having the form
\begin{equation}\label{difexp}
D=\sum_{|\alpha|+|\beta|+j+l\le m}a_{\alpha\beta j l}(r,\omega,x)
\left(-i\frac{\partial}{\partial x}\right)^\alpha
\left(-\frac{i}{r}\frac{\partial}{\partial\omega}\right)^\beta\left(-i\frac{\partial}{\partial
r}\right)^j \left(\frac{1}{r}\right)^l
\end{equation}
in a neighborhood of the edge, where $m$ is the \textit{order} of the expression and the
coefficients $a_{\alpha\beta j l}(r,\omega,x)$ are smooth functions up to $r=0$. Such
differential expressions can be realized as bounded operators
\begin{equation}\label{oper1}
D: \mathcal{W}^{s,\gamma}(\mathcal{M})\longrightarrow \mathcal{W}^{s-m,\gamma-m}(\mathcal{M})
\end{equation}
in the \emph{edge Sobolev spaces} $\mathcal{W}^{s,\gamma}$ obtained by gluing of the standard
Sobolev space $H^s$ on the smooth part of the manifold and the function space on the infinite
wedge \eqref{modell} with the norm
$$
\|u\|_{s,\gamma}=\left(\int[\xi]^{2s}\|\varkappa^{-1}_{[\xi]}\widetilde{u}\|
_{\mathcal{K}^{s,\gamma}}d\xi\right)^{1/2},\quad s,\ga\in\RR,
$$
where $[\xi]=\sqrt{1+\xi^2}$, $\widetilde{u}$ is the Fourier transform of $u$ with respect to
the variable $x$,
$$
\varkappa_\lambda u(r)=\lambda^{(k+1)/2}u(\lambda r)
$$
in a unitary action of the group $\mathbb{R}_+$ in the space $L^2(K_\Omega,d{\rm vol})$ with
the volume form $d{\rm vol}=r^k dr d\omega $ corresponding to the cone metric
$dr^2+r^2d\omega^2$ (here $k=\dim\Omega$), and the norm
$\norm{\bcdot}_{\mathcal{K}^{s,\gamma}}$ for functions on the cone $K_\Omega$ is determined by
the formula
$$
\|u\|_{\mathcal{K}^{s,\gamma}}=
\|(1+r^{-2}+\Delta_{K_\Omega})^{s/2}\rho^{s-\gamma}u\|_{L^2(K_\Omega,d{\rm vol})},
$$
in which $\Delta_{K_\Omega}$ is the Beltrami--Laplace operator with respect to the conical
metric and $\rho$ is a smooth weight function equal to $r$ in a neighborhood of $r=0$ and equal
to unity for large $r$.

\paragraph{3. Pseudodifferential operators and symbols.}

The paper~\cite{NSScS98} describes a calculus of pseudodifferential operators on manifolds with
edges extending the calculus of edge-degenerate differential operators. A pseudodifferential
operator $D$ of order $m$ in the spaces~\eqref{oper1} has a well-defined \textit{interior
symbol}, which is a function $\si(D)$ homogeneous of degree $m$ on the cotangent bundle
$T^*_0\cM$ of the manifold with edges without the zero section (the definition of
$T^*\mathcal{M}\in\Vect(M)$ can be found in the cited paper) and a well-defined \textit{edge
symbol}, which is an operator-valued function
\begin{equation}\label{resi}
    \siw(D)(x,\xi):\cK^{s,\ga}(K_\Om)\lra \cK^{s-m,\ga-m}(K_\Om),\quad (x,\xi)\in T_0^*X,
\end{equation}
in the function space on the infinite wedge. The edge symbol possesses the \textit{twisted
homogeneity} property
\begin{equation}\label{scrudge}
    \siw(D)(x,\la\xi)=\la^m\ka_\la\siw(D)(x,\xi)\ka_\la^{-1},\quad \la\in\RR_+.
\end{equation}
In particular, the edge symbol of the operator~\eqref{difexp} is equal to
\begin{equation*}
    D=\sum_{|\alpha|+|\beta|+j+l= m}a_{\alpha\beta j l}(0,\omega,x)
\xi^\alpha
\left(-\frac{i}{r}\frac{\partial}{\partial\omega}\right)^\beta\left(-i\frac{\partial}{\partial
r}\right)^j \left(\frac{1}{r}\right)^l.
\end{equation*}
The notion of interior and edge symbols is important in that the following assertion holds.
\begin{proposition}\label{kritkomp}
A pseudodifferential operator $D$ of order $m$ in the spaces~\eqref{oper1} is compact if and
only if $\si(D)=0$ and $\siw(D)=0$.
\end{proposition}

\paragraph{4. The composition theorem and ellipticity.}
The main property of the calculus of pseudodifferential operators on manifolds with edges is
expressed by the composition theorem.
\begin{theorem}\label{product-gen}
The composition of edge-degenerate pseudodifferential operators corresponds to the composition
of their symbols \rom(interior and edge\rom{):}
$$
\sigma(D_1D_2)=\sigma(D_1)\sigma(D_2),\quad \siw(D_1D_2)=\siw(D_1)\siw(D_2).
$$
\end{theorem}
In conjunction with the compactness criterion given by proposition~\ref{kritkomp}, the
composition formula results in the following finiteness theorem.
\begin{theorem}\label{finita}
If a pseudodifferential operator $D$ of order $m$ acting in the space~\eqref{oper1} is elliptic
\rom(i.e., its interior symbol is invertible everywhere outside the zero section on $T^*\cM$
and the edge symbol is invertible in the spaces~\eqref{resi} everywhere outside the zero
section on $T^*X$\rom), then it is Fredholm.
\end{theorem}

\medskip

Naturally, all the preceding is also valid for operators acting in spaces of sections of vector
bundles on $M$.

\paragraph{5. Order reduction.}

In this paper, we are interested in the classification of elliptic pseudodifferential
operators~\eqref{oper1} modulo stable homotopies. (The precise definition of stable homotopy
will be given below.) Let us show that this problem can actually be reduced to the special case
of zero-order operators in the spaces $\cW^{s,\ga}$ for $s=\ga=0$. Indeed, an operator $D$ is
elliptic for given $s$ (and fixed $\ga$) if and only if it is elliptic for any $s$ (with the
same $\ga$); see, e.g., \cite{NSScS98}. Consequently, a homotopy in the class of elliptic
operators for some $s$ is valid for all $s$. Hence without loss of generality we can assume
that $s=\ga$, i.e., consider the operators
\begin{equation}\label{oper2}
D: \mathcal{W}^{\ga,\gamma}(\mathcal{M})\longrightarrow
\mathcal{W}^{\ga-m,\gamma-m}(\mathcal{M}).
\end{equation}
Next, there exist elliptic operators
\begin{equation}\label{redu}
    V:\cW^{0,0}(\cM)\lra \cW^{\ga,\ga}(\cM),\quad
    \wt V:\cW^{\ga-m,\ga-m}(\cM)\lra \cW^{0,0}(\cM)
\end{equation}
of index zero. Then the mapping
\begin{equation*}
    D\longmapsto \wt V D V
\end{equation*}
reduces elliptic operators (and homotopies) in the spaces~\eqref{oper2} to those in the space
$\cW^{0,0}(\cM)$. The inverse mapping (modulo compact operators, which plays no role if
homotopies in the class of elliptic operators are considered) naturally has the form
\begin{equation*}
    D\longmapsto \wt V^{-1} D V^{-1},
\end{equation*}
where $\wt V^{-1}$ and $V^{-1}$ are almost inverses of $\wt V$ and $V$, respectively.

\paragraph{6. Pseudodifferential operators of order zero.}
By virtue of order reduction, in what follows we are mainly interested in pseudodifferential
operators of order zero in the space $\cW^{0,0}(\cM)$. Hence we give a more detailed
description of their construction and properties, mainly following~\cite{NSScS98} with some
simplifications (related to the fact that the paper~\cite{NSScS98} deals with operators of
arbitrary order and not only compactness, but also smoothing properties of remainders in
composition formulas are taken into account). For simplicity, we present all facts for
operators acting in function spaces on $\cM$. The generalization to operators acting in spaces
of sections of vector bundles on $\cM$ is trivial.

\subparagraph{Edge symbols.} First, we describe the class of edge symbols used here.
\begin{definition}\label{resid}
An \textit{edge symbol} is a family $D(x,\xi)$, $(x,\xi)\in T_0^*X$, of operators in function
spaces on the cones $K_\Om$ with the following properties.
\begin{enumerate}
    \item For any multi-indices $\a,\be$, $\abs{\a}+\abs{\be}=0,1,2,\dots$\,,
          the derivatives
\begin{equation*}
    D^{(\a,\be)}(x,\xi):\cK^{0,0}(K_\Om)\lra\cK^{0,0}(K_\Om)
\end{equation*}
    are continuous operators.
    \item The \textit{twisted homogeneity condition}
\begin{equation*}
    D(x,\la\xi)=\ka_\la D(x,\xi)\ka_\la^{-1}, \quad\la\in\RR_+,
\end{equation*}
          holds.
    \item Modulo compact operators, one has the representation
\begin{equation}\label{votetoda}
    D(x,\xi)=d\BL(x,\xi\ovs2r,i\ovs1{r\pd{}r}\BR),
\end{equation}
    where $d(x,\eta,p)$ is a classical pseudodifferential operator with parameters
    $(\eta,p)\in T^*_xX\times\cL_{-(k+1)/2}$ in the sense of Agranovich--Vishik~\cite{AgVi1}
    of order zero on $\Om$ smoothly depending on the additional parameter
    $x\in X$.
\end{enumerate}
\end{definition}
Here $\cL_{-(k+1)/2}=\{\im p=-(k+1)/2\}$ is the weight line, and the function of the operator
$i{r\pa/\pa r}$ in~\eqref{votetoda} is defined with the help of the Mellin transform on this
weight line.

\begin{definition}\label{vot-vam}
The \textit{interior symbol} of the edge symbol $D(x,\xi)$ is the principal symbol
$\si(D)=\si(d)$ in the sense of Agranovich--Vishik of the corresponding pseudodifferential
operator $d(x,\eta,p)$. The \textit{conormal symbol} of $D(x,\xi)$ is the operator family
$\si_c(D)=d(x,0,p)$ in the space $L^2(\Om)$.
\end{definition}

The main properties of edge symbols are expressed by the following theorem.

\begin{theorem}
The following assertions hold.
\begin{enumerate}
    \item[\rom{1.}] Definition~\rom{\ref{resid}} is consistent\rom; i.e., the
                    operator~\eqref{votetoda}
                    is always bounded in $\cK^{0,0}(K_\Om)$.
    \item[\rom{2.}] The operator~\eqref{votetoda} is compact if and only if its
    interior and conormal symbol are zero. \rom(In particular, it follows that the
    interior and conormal symbol of an edge symbol are well defined.\rom)
    \item[\rom{3.}] Edge symbols form a local $C^*$-algebra, and, modulo compact edge symbols,
    the product of pseudodifferential operators $d(x,\eta,p)$ corresponds to the product of
    the respective edge symbols $D(x,\xi)$ and the adjoint operator corresponds to the adjoint
    edge symbol.
    \item[\rom{4.}] \rom(Corollary.\rom) The mapping that takes each edge symbol to its
    interior and conormal symbols is linear and multiplicative and commutes with the passage
    to the adjoint operator.
    \item[\rom{5.}] \rom(Norm estimates modulo compact operators.\rom)
    For an edge symbol $D(x,\xi)$ of order zero, one has\footnote{Here and in the following,
    $\mathcal{K}$ is either the ideal of compact operators or (where there is a bundle) a bundle
    of algebras of compact operators acting in function spaces on the fibers.
    Likewise, $\mathcal{B}$ corresponds to the algebra of bounded operators.}
\begin{multline}
\label{est2} \inf_{K\in C(S^*X,\mathcal{K})}\max_{S^*X}\|D+K\|_{\mathcal{B}
(\cK^{0,0}(K_\Omega))} =\max\left(\max_{\partial S^*\mathcal{M}} |\sigma(D)|,\sup_{X\times
\mathbb{R}} \|\sigma_c(D)\|_{\mathcal{B}(L^2(\Omega))}\right).
\end{multline}
    \item[\rom{6.}] The commutator $[D(x,\xi),\ph]$ is compact for any continuous function
                    $\ph(r)$ on $\RR_+$ equal to zero for sufficiently large $r$.
    \item[\rom{7.}] The product $\ph(r)D(x,\xi)$ is a symbol of order zero
                    and has a compact fiber variation \rom{\cite{Luk1}} on
                    $T^*_0X$ if $\ph(r)$ is the same as in item~\rom6.
    \item[\rom{8.}] For any pair \rom(interior symbol, conormal symbol\rom) satisfying
    the compatibility condition \rom($\sigma(d)$ for $\eta=0$ is equal to
    $\sigma(\sigma_c)$\rom) one can construct an edge symbol.
\end{enumerate}
\end{theorem}

The algebra of zero-order edge symbols will be denoted by $\Psi_\wedge(X)$.

\subparagraph{Pseudodifferential operators.}

Now we can describe the class of zero-order pseudodifferential operators.

Let $\cA$ be the algebra of classical zero-order pseudodifferential operator $A$ on the open
manifold $\ciM$ with the following properties:
\begin{enumerate}
    \item[1)] the principal symbol $\si(A)$ is a smooth function on $T^*_0\cM$
              up to the boundary;
    \item[2)] the operator $A$ is continuous in the space $\cW^{0,0}(\cM)$;
    \item[3)] $A$ compactly commutes with $C(\cM)$.
\end{enumerate}

\begin{definition}
A \textit{zero-order pseudodifferential operator} on $\cM$ is a continuous operator
\begin{equation}\label{psiddo}
    B:\cW^{0,0}(\cM)\lra \cW^{0,0}(\cM)
\end{equation}
representable modulo compact operators in the form
\begin{equation}\label{mainformula}
    B=(\ph(r)D)\BL(\ovs2x,-i\ovs1{\pd{}x}\BR)+rA,
\end{equation}
where $D(x,\xi)\in\Psi_\wedge(X)$ is an edge symbol, $\varphi(r)$ is a smooth function on $M$
equal to $1$ on $\pa M$ and zero outside a sufficiently small neighborhood of the boundary, and
$A\in\cA$ is a classical pseudodifferential operator on $\ciM$. (Here $r$ is the distance to
$\pa M$.)
\end{definition}

\begin{definition}
The \textit{edge symbol} of the operator~\eqref{psiddo} is the operator family
\begin{equation}\label{resi2}
    \siw(B)=D(x,\xi).
\end{equation}
The \textit{interior symbol} of the operator~\eqref{psiddo} is the function
\begin{equation}\label{vnusi}
    \si(B)=\ph(r)\si(D)+r\si(A).
\end{equation}
\end{definition}

The set of zero-order pseudodifferential operator on $\cM$ will be denoted by $\Psi(\cM)$. The
main properties of pseudodifferential operators are expressed by the following theorem.

\begin{theorem}\label{mainprop}
The following assertions hold.
\begin{enumerate}
    \item[\rom{1.}] The definition is consistent\rom; i.e., the operator~\eqref{mainformula}
                    is always bounded in $\cW^{0,0}(\cM)$.
    \item[\rom{2.}] The interior and edge symbols of the operator~\eqref{mainformula}
                    are well defined. It is compact if and only if both symbols are zero.
    \item[\rom{3.}] The interior and edge symbols satisfy the compatibility condition
\begin{equation}\label{soglass}
    \si(B)|_{\pa T^*_0\cM}=\si(\siw(B)).
\end{equation}
    For any pair \rom(interior symbol, edge symbol\rom) satisfying the compatibility condition,
    one can construct the corresponding pseudodifferential operator.
    \item[\rom{4.}] Pseudodifferential operators form a local $C^*$-algebra,
    and the mapping taking each pseudodifferential operator to its
    interior and edge symbols is a $*$-homomorphism.
    \item[\rom{5.}] \rom(Norm estimates modulo compact operators.\rom) For a zero-order
    edge-degenerate pseudodifferential operator $D$ one has
\begin{equation}
\inf_{K\in \mathcal{K}}\|D+K\|_{\mathcal{B}( \cW^{0,0}(\mathcal{M})
)}=\max\left(\max_{S^*\mathcal{M}} |\sigma(D)|,\max_{S^*X}
\|\sigma_\Lambda(D)\|_{\mathcal{B}(\cK^{0,0}(K_\Omega))}\right)
\end{equation}
    \item[\rom{6.}] If $B\in\Psi(\cM)$, then the commutator $[B,\ph]$
                    is compact for any function $\ph\in C(\cM)$.
\end{enumerate}
\end{theorem}

\subparagraph{The norm closure.}

The norm estimates modulo compact operators imply the following description of the Calkin
algebras of the closures $\overline{\Psi(\mathcal{M})}$ and $\overline{\Psi_\Lambda(X)}$.
\begin{corollary}
The interior and edge symbol homomorphisms \rom(for operators on $\mathcal{M}$\rom)
$$
\Psi(\mathcal{M})\stackrel{\sigma}\to C^\infty(S^*\mathcal{M}),\;\;
\Psi(\mathcal{M})\stackrel{\sigma_\Lambda}\to \Psi_\Lambda(X),
$$
induce the isomorphism
$$
\overline{\Psi(\mathcal{M})}/\mathcal{K}\simeq\left\{ (a,a_\Lambda) \Bigl|
\begin{array}{c}
a\in C(S^*\mathcal{M}),\quad a_\Lambda\in \overline{\Psi_\Lambda(X)}:\vspace{1mm}\\
a|_{\partial S^*\mathcal{M}}=\sigma(a_\Lambda)
\end{array}
\right\}.
$$
The interior and conormal symbol homomorphisms \rom(for edge symbols\rom)
$$
\Psi_\Lambda (X)\stackrel{\sigma}\to C^\infty(\partial S^*\mathcal{M}),\;\; \Psi_\Lambda
(X)\stackrel{\sigma_c}\to \Psi_c(X),
$$
where $\Psi_c$ is the algebra of conormal symbols, induce the isomorphism
$$
\overline{\Psi_\Lambda(X)}/C( S^*X,\mathcal{K})\simeq \left\{ (a,a_c) \Bigl|
\begin{array}{c}
a\in C(\partial S^*\mathcal{M}),\quad a_c\in \overline{\Psi_c(X)}:\vspace{1mm}\\
a|_{S(T^*X\oplus \mathbf{1})}=\sigma(a_c)
\end{array}
\right\}.
$$
\end{corollary}

\section{An Element in the $K$-Homology of the Singular Space}
In this section, we show how an elliptic operator on a manifold with edges gives rise to an
element in the analytic $K$-homology of the space $\mathcal{M}$. (A detailed exposition of the
theory of analytic $K$-homology can be found in \cite{HiRo1}, \cite{Bla1}, and \cite{Kas3}. An
introduction to the theory can be found in \cite{SaSt15}.) To the best of the authors'
knowledge, this correspondence was used for the first time in~\cite{Has1}. In accordance with
the preceding, we consider only zero-order operators.

Let
$$
D:\mathcal{W}^{0,0}(\mathcal{M},E)\longrightarrow \mathcal{W}^{0,0}(\mathcal{M},F)
$$
be an elliptic operator of order zero in sections of vector bundles $E,F\in\Vect(M)$. The
commutator $[D,f]$ with a function $f\in C^\infty(M)$ is compact if the restrictions of the
function to the fibers of $\pi$ are constant functions. (This follows from the composition
formula.) Thus $D$ is a generalized elliptic operator on $\mathcal{M}$ in the sense of Atiyah
\cite{Ati4} and hence defines a class in the analytic $K$-homology $K_0(\mathcal{M})$ of the
singular space $\cM$. Let us give a precise construction of the corresponding cycle.

If $D$ is self-adjoint (and $E=F$), then we consider the \textit{normalization}
\begin{equation}\label{oda1}
\mathcal{D}=(P_{\ker D}+D^2 )^{-1/2}D:\mathcal{W}^{0,0}(\mathcal{M},E)\longrightarrow
\mathcal{W}^{0,0}(\mathcal{M},E),
\end{equation}
where $P_{\ker D}$ is the orthogonal projection on the kernel.

In the general case, we consider the self-adjoint operator
\begin{equation}\label{eva1}
\mathcal{D}=\left(
\begin{array}{cc}
0 &  \!\!\!\!\!\! D(P_{\ker D}+D^*D)^{-1/2} \\
(P_{\ker D}+D^*D)^{-1/2}D^*\!\!\!\!\!\!& 0
\end{array}
\right) :\mathcal{W}^{0,0}(\mathcal{M},E\oplus F)\rightarrow
\mathcal{W}^{0,0}(\mathcal{M},E\oplus F),
\end{equation}
which is odd with respect to the $\mathbb{Z}_2$-grading of the space
$\mathcal{W}^{0,0}(\mathcal{M},E)\oplus \mathcal{W}^{0,0}(\mathcal{M},F)$. By $C(\mathcal{M})$
we denote the algebra of continuous functions on $\mathcal{M}$.

\begin{proposition}\label{modul}
The operators \eqref{oda1} and \eqref{eva1} are zero-order elliptic pseudodifferential
operators and define elements in $K$-homology\rom; these elements will be denoted by
$$
[D]\in K_*(\mathcal{M}),
$$
where $*=1$ for self-adjoint operators and $*=0$ in the general case.
\end{proposition}

\begin{proof}
1. The operator $D^*D$ is pseudodifferential, and the same is true for $P_{\ker D}$, since the
latter operator is finite rank and hence compact. Thus $D^*D+P_{\ker D}$ is a
pseudodifferential operator, and since it is invertible, it follows that the inverse is also a
pseudodifferential operator (recall that $\ov{\Psi(\cM)}$ is a $C^*$-algebra). To prove that
$(D^*D+P_{\ker D})^{-1/2}$ is a pseudodifferential operator, it remains to use the formula
$$
A^{-1/2}=\frac 1\pi \int_0^\infty \lambda^{-1/2}(A+\lambda)^{-1}d\lambda
$$
for a self-adjoint strongly positive operator $A$ (e.g., see~\cite[p.~165]{Bla1}).

2. The operators $\mathcal{D}$  in \eqref{oda1} and \eqref{eva1} are self-adjoint operators
acting in a $*$-module over the $C^*$-algebra $C(\mathcal{M})$  and have the properties
$$
[\mathcal{D},f]\in \mathcal{K},\qquad (\mathcal{D}^2-1)f\in \mathcal{K}.
$$
Thus we have the Fredholm modules
\begin{equation}
\label{modules} [D]=\left\{
\begin{array}{rl}
[\mathcal{D},\mathcal{W}^{0,0}(\mathcal{M},E)]\in K_1(\mathcal{M})
& \text{ if } D=D^*\vspace{2mm}\\
\left[\mathcal{D},\mathcal{W}^{0,0}(\mathcal{M},E\oplus F)\right]\in K_0(\mathcal{M}) &
\text{ if } D\ne D^*.\\
\end{array}
\right.
\end{equation}
\end{proof}

\begin{remark}
If only the interior symbol of $D$ is elliptic,  then there is a well-defined element in the
$K$-homology of the open smooth part of $\cM$:
\begin{equation}
[D]\in K_*(\mathcal{M}\setminus X).
\end{equation}
Here in the definition of the elements \eqref{oda1} and \eqref{eva1} one should replace the
expression $(P_{\ker D}+D^*D)^{-1/2}$ by an arbitrary self-adjoint edge-degenerate
pseudodifferential operator with interior symbol $(\sigma(D)^*\sigma(D))^{-1/2}$.
\end{remark}

\section{The Homotopy Classification}

We recall the standard equivalence relation on the set of pseudodifferential operator acting in
sections of vector bundles, namely, \textit{stable homotopy}.

\begin{definition}
Two operators
$$
D:\mathcal{W}^{0,0}(\mathcal{M},E)\to \mathcal{W}^{0,0}(\mathcal{M},F),\quad
D':\mathcal{W}^{0,0}(\mathcal{M},E')\to \mathcal{W}^{0,0}(\mathcal{M},F')
$$
are said to be \textit{stably homotopic} if they are homotopic modulo stabilization by vector
bundle isomorphism, i.e., if there exists an continuous homotopy of elliptic operators
$$
D\oplus 1_{E_0}\sim f^*\bigl(D'\oplus 1_{F_0}\bigr)e^*,
$$
where $E_0,F_0\in\Vect(M)$  are vector bundles and
$$
e:E\oplus E_0\longrightarrow E'\oplus F_0,\qquad f:F'\oplus F_0\longrightarrow F\oplus E_0
$$
are vector bundle isomorphisms.
\end{definition}

Stable homotopy is an equivalence relation on the set of all elliptic edge-degenerate
pseudodifferential operators acting in sections of vector bundles. By $\Ell_0(\mathcal{M})$ we
denote the set of elliptic operators modulo stable homotopies. This set is a group with respect
to the direct sum of elliptic operators. The inverse element corresponds to an almost inverse
(i.e., an inverse modulo compact operators), and the unit is the equivalence class of trivial
operators.

In a similar way, one defines odd elliptic theory $\Ell_1(\mathcal{M})$ as the group of stable
homotopy equivalence classes of elliptic self-adjoint operators. Here the class of trivial
operators consists of Hermitian isomorphisms of vector bundles.

\medskip

The \textit{homotopy classification problem for elliptic operators} on the manifold
$\mathcal{M}$ is the problem of computing the groups $\Ell_*(\mathcal{M})$.

\medskip

The following theorem solves the classification problem for manifolds with edges and is the
main result of the paper.
\begin{theorem}\label{thmain}
There is an isomorphism
$$
\Ell_*(\mathcal{M}) \stackrel{\varphi}\simeq K_*(\mathcal{M}),
$$
which takes each elliptic operator $D$ to the element defined in Proposition~\rom{\ref{modul}}.
\end{theorem}

\begin{corollary}
Two elliptic operators $D_1$ and $D_2$ are stably rationally homotopic if and only if they have
the same indices with coefficients in an arbitrary bundle on $\mathcal{M}$:
\begin{equation}\label{indpair} \ind
(1\otimes p)(D_1\otimes 1_N)(1\otimes p)=\ind (1\otimes p)(D_2\otimes 1_N)(1\otimes p),
\end{equation}
where $p\in \Mat(N\times N,C(\mathcal{M}))$ is a matrix projection.
\end{corollary}
This follows from the nondegeneracy (on the free parts of the groups) of the natural pairing
$K_0(\mathcal{M})\times K^0(\mathcal{M})\longrightarrow \mathbb{Z}$, which is just defined by
the formula~\eqref{indpair}.

\begin{proof}[Proof of the theorem]
The mapping is well defined, since homotopies of elliptic operators give rise to homotopies of
the corresponding Fredholm modules, i.e., result in the same element in $K$-homology. Bundle
isomorphisms give degenerate modules.

\smallskip

Let us now prove that the mapping is an isomorphism. We split the proof into three stages.

\paragraph{1. Reduction of $\Ell$-groups to $K$-groups of $C^*$-algebras (see~\cite{Sav8}).}

We interpret edge-degenerate operators in sections of vector bundles as operator generated by
the pair of algebras
$$
C^\infty(M)\subset \Psi(\mathcal{M})
$$
of scalar operators. The embedding corresponds to the conventional action of functions as
multiplication operators. Namely, an arbitrary edge-degenerate pseudodifferential operator of
order zero can be represented in the form
$$
D':\im P\longrightarrow \im Q,
$$
where $P=P^2$ and $Q=Q^2$ are matrix projections with coefficients in the function algebra
$C^\infty(M)$ and $D'$ is a matrix operator whose entries belong to the operator algebra
$\Psi(\mathcal{M})$.

We obtain a group isomorphic to $\Ell(\mathcal{M})$ if, instead of operators with smooth
symbols, we consider operators whose symbols are only \emph{continuous}, i.e., pass to the
closure $\overline{\Psi(\mathcal{M})}$ of the algebra of pseudodifferential operators with
respect to the operator norm. (The fact that these groups are isomorphic follows from
Theorem~\ref{mainprop}). By $\Sigma\stackrel{\rm def}=\overline{\Psi(\mathcal{M})}/\mathcal{K}$
we denote the algebra of continuous symbols.

The results of~\cite{Sav8} give the isomorphisms\footnote{Note that if the edges are absent,
then the isomorphism is just the Atiyah--Singer difference construction \cite{AtSi1}.}
$$
\Ell_*(\mathcal{M})\stackrel\chi\simeq K_*(\Con(C(M)\to \Sigma)).
$$
Here
$$
\Con(A\stackrel{f}\longrightarrow B)=\Bigl\{(a,b(t))\in A\oplus C_0([0,1),B)\;|\;
f(a)=b(0)\Bigr\}
$$
is the cone of the algebra homomorphism $f:A\to B$. In the odd case, one can rewrite the
$K$-group in the form
$$
K_1(\Con(C(M)\to \Sigma))\simeq K_{0}(\Sigma)/K^{0}(M).
$$
The composition of the last isomorphism with $\chi$ is a generalization of the
Atiyah--Patodi--Singer isomorphism \cite{APS3}; i.e., self-adjoint elliptic operators modulo
stable homotopy are isomorphic to symbols-projections modulo projections determining sections
of bundles.

\paragraph{2. A diagram relating $K$-theory of algebras and $K$-homology.}

Consider the commutative diagram
\begin{equation}
\label{dia1}
\begin{array}{rcccccl}
0\rightarrow & \Sigma_0 & \longrightarrow & \Sigma & \stackrel\sigma \longrightarrow &
C(S^*\mathcal{M}) &
\to 0\\
& \uparrow & & \uparrow & & \uparrow \\
0\rightarrow & 0 & \longrightarrow & C(M) & = & C(M) & \to 0
\end{array}
\end{equation}
with exact rows. Here $\sigma$ is the interior symbol and $\Sigma_0=\ker \sigma$ is the
corresponding ideal. The diagram induces the exact sequence
\begin{equation}\label{conu1}
0\to S\Sigma_0 \longrightarrow \Con(C(M)\to \Sigma) \longrightarrow  C_0(T^*\mathcal{M}) \to 0
\end{equation}
of the cones of the vertical homomorphisms. Here $S$ stands for the suspension:
$S\Sigma_0=C_0((0,1),\Sigma_0)$.

For brevity, we set ${A}=\Con(C(M)\to \Sigma)$.

The key point of the proof is the construction of the commutative diagram
\begin{equation}
\label{dia3}
\begin{array}{ccccccccc}
K^{*+1}_c(T^*\mathcal{M}) & \stackrel\partial \to & K_*(S\Sigma_0) &\to & K_*({A}) &\to &
K^*_c(T^*\mathcal{M}) &
\stackrel\partial\to & K_{*+1}(S\Sigma_0) \\
\downarrow & \boxed{A} & \downarrow & \boxed{B} & \downarrow & \boxed{C} & \downarrow &
\boxed{D} & \downarrow\\
K_{*+1}(M\setminus \partial M) &\stackrel\partial \to & K_*(X) & \to & K_*(\mathcal{M}) & \to &
K_*(M\setminus \partial M) & \to & K_{*+1}(X)
\end{array}
\end{equation}
relating the exact sequence induced by \eqref{conu1} in $K$-theory of algebras to the exact
sequence of the pair $X\subset \mathcal{M}$ in $K$-homology.

The vertical arrows in \eqref{dia3} are induced by quantizations. Namely, the elements of
$K$-groups in the upper row correspond to some symbols, and the vertical mappings take these
symbols to the corresponding operators. More precisely,
\begin{itemize}
\item the mappings $K^*_c(T^*\mathcal{M})\to K_*(M\setminus \partial M)$ and $K_*(A)\to
K_*(\mathcal{M})$ are determined by quantizations (see the preceding section);

\item the mappings $K_*(S\Sigma_0)\longrightarrow K_*(X)$ are induced by quantization of edge
symbols. In more detail, the mapping $K_1(\Sigma_0)\longrightarrow K_{0}(X)$ takes an edge
symbol $1+u(x,\xi)$, $u\in\Sigma_0$, invertible on $S^*X$ to the operator
\begin{equation}\label{luke}
1+u\left(x,-i\frac\partial{\partial x}\right)
\end{equation}
in the space $L^2(X,\mathcal{K}^{0,0}(K_\Omega))$ on the infinite wedge. The quantization is
well defined, since $u(x,\xi)$ has a compact fiber variation (the interior symbol of the edge
symbol is zero). The mapping $K_0(\Sigma_0)\longrightarrow K_{1}(X)$ takes a self-adjoint edge
symbol-projection $p(x,\xi)$ to the operator
\begin{equation}\label{oda2}
2p\left(x,-i\frac\partial{\partial x}\right)-1.
\end{equation}
\end{itemize}

\paragraph{3. An application of the five lemma.}
Once we prove that the diagram \eqref{dia3} commutes, it follows from the five lemma that the
middle vertical arrow is an isomorphism, since the quantizations $K^*_c(T^*\mathcal{M})\to
K_*(M\setminus\partial M)$ are defined on the interior of the smooth manifold $M$ with boundary
and are isomorphisms (e.g., see~\cite{Kas3}). The isomorphism of edge quantizations
$K_*(S\Sigma_0)\to K_*(X)$ will be established below in Sec.~\ref{tri}.

The fact that the diagram commutes can be established as follows:
\begin{itemize}
\item The square $\boxed{C}$ in \eqref{dia3} commutes automatically, since the horizontal
arrows are forgetful mappings;
\item The fact that the square $\boxed{B}$ commutes will be proved in Sec.~\ref{tri}.
\item The fact that the square $\boxed{A}$, including the boundary maps, commutes will be
verified in Sec.~\ref{chetire}.
\end{itemize}
This completes the proof of the theorem up to the above-mentioned computations.
\end{proof}

\section{Computations on the Edge}\label{tri}

\paragraph{1. The isomorphism $K_*(\Sigma_0)\simeq K_{*+1}(X)$.}
Consider the diagram
\begin{equation}
\label{dia2}
\begin{array}{ccccccc}
\rightarrow & K^*(S^*X) & \rightarrow & K_*(\Sigma_0)  &
 \rightarrow  & K^*_c(X\times \mathbb{R}) & \rightarrow\\
& \parallel &  & \downarrow {L} & & \parallel \\
\rightarrow & K^*(S^*X) & \rightarrow & K^{*+1}_c(T^*X)&
 \rightarrow  & K^*_c(X\times \mathbb{R}) & \rightarrow,
\end{array}
\end{equation}
which compares the $K$-theory sequence corresponding to the short exact sequence
\begin{equation}
\label{solv2} 0\to C(S^*X,\mathcal{K})\longrightarrow \Sigma_0
\stackrel{\sigma_c}\longrightarrow C_0(X\times \mathbb{R},\mathcal{K}) \to 0
\end{equation}
(where $\sigma_c$ is the conormal symbol) with the sequence of $K$-groups of the pair
$S^*X\subset B^*X$ formed by the unit sphere and ball bundles in $T^*X$. The mapping $L$ is the
difference constructions for pseudodifferential operators \eqref{luke}, \eqref{oda2} with
operator-valued symbols in the sense of Luke (see~\cite{Luk1} and~\cite{NSScS10}).

\begin{lemma}
The diagram \eqref{dia2} commutes.
\end{lemma}

\begin{proof}
1. The commutativity of the squares
$$
\begin{array}{ccc}
K^*(S^*X) & \longrightarrow &  K_*(\Sigma_0) \\
\parallel & & \downarrow {L}\\
K^*(S^*X) & \longrightarrow &  K^{*+1}(T^*X)
\end{array}
$$
follows from the fact that for finite-dimensional symbols the difference constriction coincides
with the Atiyah--Singer difference constriction.

2. The commutativity of the squares
$$
\begin{array}{ccc}
K_*(\Sigma_0) & \stackrel{\sigma_c}\longrightarrow & K^*_c(X\times \mathbb{R})  \\
\downarrow {L} & & \parallel\\
K^{*+1}_c(T^*X)& \stackrel{j^*}\longrightarrow & K^{*+1}(X),
\end{array}\quad j:X\to T^*X,
$$
follows from the index formula \cite{NSScS3}
\begin{equation}\label{indf1}
\beta \ind D_y=\ind \sigma_c(D_y)\in K^1_c(Y\times \mathbb{R})
\end{equation}
for a family of elliptic operators $D_y,$ $y\in Y$, with unit interior symbol on the infinite
cone. Here $Y$ is a compact parameter space and $\beta$ is the periodicity isomorphism
$K(Y)\simeq K^1_c(Y\times \mathbb{R})$.

The formula \eqref{indf1} applies directly to the group $K_1(\Sigma_0)$, and for the
$K_0$-group one uses the suspension (cf.~\eqref{oda2}).

3. The commutativity of the squares
$$
\begin{array}{ccc}
 K^{*+1}_c(X\times \mathbb{R}) & \stackrel\partial \longrightarrow & K^*(S^*X)  \\
\parallel & & \parallel\\
K^*(X)& \stackrel{p^*}\longrightarrow &  K^*(S^*X)
\end{array}, \quad p:S^*X\to X,
$$
also follows from the above-mentioned index formula, since the boundary mapping in $K$-theory
of algebras is an index mapping. We leave details to the reader.
\end{proof}

By applying the five lemma, we arrive at the desired corollary.
\begin{corollary}
The quantization $K_*(\Sigma_0)\to K_{*+1}(X)$ is an isomorphism\rom: $K_*(\Sigma_0)\simeq
K^{*+1}_c(T^*X)\simeq K_{*+1}(X)$.
\end{corollary}

\paragraph{2. Commutativity of the square $\boxed{B}$.}
To be definite, we consider the even case. The odd case can be treated in a similar way.

The image of the composite mapping $K_0(S\Sigma_0)\to K_0({A})\to K_0(\mathcal{M})$ corresponds
to elliptic operators on $\mathcal{M}$ of the form $1+\mathbf{G}$, where $\mathbf{G}$ is an
operator with zero interior symbol. We must show that the element
$$
[1+\mathbf{G}:\mathcal{W}^{0,0}(\mathcal{M})\to \mathcal{W}^{0,0}(\mathcal{M})]\in
K_0(\mathcal{M})
$$
coincides with the element
$$
\left[1+g\left(x,-i\frac\partial{\partial x} \right): \mathcal{W}^{0,0}({W})\longrightarrow
\mathcal{W}^{0,0}({W})\right] \in K_0(\mathcal{M})
$$
determined by the operator on the infinite wedge $W$ with symbol
$g(x,\xi)=\sigma_\Lambda(\mathbf{G})$. For the latter operator, the module structure on the
spaces is defined as follows: a function $f\in C(\mathcal{M})$ acts as the multiplication by
its restriction to the edge.

The equality of these two elements can be established in two steps:
\begin{enumerate}
\item (cutting away the smooth part of the manifold) without changing the element of the
$K$-homology group, one can proceed to the restriction of the operator $1+\mathbf{G}$ to a
neighborhood of the edge (the original operator and its restriction define stably equivalent
Fredholm modules and hence the same elements in $K_0(\mathcal{M})$);

\item (homotopy of the module structure) in the neighborhood of the edge, the module structure
$f,u\mapsto fu$ cam be homotopied to  $f,u\mapsto f(x,0)u(x,\omega,r)$ by scaling represented
by the formula $f(x,\omega,\varepsilon r)u(x,\omega,r)$ in local coordinates $(x,\omega,r)$.
\end{enumerate}
We have established that the square \fbox{B} commutes.

\section{Comparison of the Boundary Mappings}\label{chetire}

Let us verify that the boundary mappings in $K$-theory of algebras and $K$-homology of spaces
are compatible.

\paragraph{1.}
First, we show that the boundary mappings can be defined in terms of the restriction of
structures to the boundary. Consider the commutative diagram of $C^*$-algebras
$$
\begin{array}{rcccccl}
0\to & S\Sigma_0 & \to & {A} &\to  & C_0(T^*\mathcal{M})& \to 0\\
& \parallel & & \downarrow & & \downarrow & \\
0\to & S\Sigma_0 & \to &  \Con(C(\partial M)\to \overline{\Psi_\Lambda(X)})  &\to  & C_0(\partial
T^*\mathcal{M})& \to 0
\end{array}
$$
relating two short exact sequences.  Next, we consider the diagram of spaces and continuous
mappings
$$
\begin{array}{ccc}
\partial M & \subset & M \\
\pi\downarrow\;\;\; & & \downarrow \\
X & \subset & \mathcal{M}.
\end{array}
$$
Since the boundary mapping is natural, we have the following lemma.
\begin{lemma}
The diagrams
\begin{align}
&\begin{array}{ccc}
K^*_c(T^*\mathcal{M}) & \stackrel\partial\longrightarrow & K_{*+1}(S\Sigma_0)\\
\downarrow & &\parallel\\
K^*_c(\partial T^*\mathcal{M}) & \stackrel{\partial'}\longrightarrow & K_{*+1}(S\Sigma_0),
\end{array}\\[2mm]
&\begin{array}{ccc}
K_*(M\setminus\partial M) & \stackrel\partial\longrightarrow& K_{*+1}(X)\\
\partial'\downarrow & &\parallel\\
K_{*+1}(\partial M) & \stackrel{\pi_*}\longrightarrow & K_{*+1}(X)
\end{array},\quad\pi:\partial M\longrightarrow X,
\end{align}
commute.
\end{lemma}

\paragraph{2.}
To compute the boundary mappings $K^*_c(\partial T^*\mathcal{M})
\stackrel{\partial'}\longrightarrow K_{*+1}(S\Sigma_0)$, we reduce the algebra of edge symbols
to the simpler algebra of operator families on the fibers $\Omega$ with parameters in the sense
of Agranovich--Vishik \cite{AgVi1}). We relate these two algebras by the semiclassical
quantization method (e.g., see~\cite{NScS12}).

In our case, semiclassical quantization is a family of linear mappings (see the Appendix)
$$
T_h:\Psi_{T^*X\times \mathbb{R}}(\Omega) \longrightarrow \Sigma|_{\partial M}, \quad h\in (0,1],
$$
\begin{equation}\label{quasi}
 (T_h u)(\xi)=u\BL(\stackrel{2} r\xi,ih \stackrel 1
{r\frac\partial{\partial r}}+ih\frac{n+1}2\BR),\quad \xi\in T^*X,
\end{equation}
where by $\Psi_{T^*X\times \mathbb{R}}(\Omega) $ we denote the algebra of smooth families of
pseudodifferential operators on the fibers $\Omega$ with parameters in $T^*X\times \mathbb{R}$;
these mappings satisfy the following relations as $ h\to 0$:
$$
T_h (ab)= T_h(a)T_h(b)+o(1),\;\;\; \left(T_h(a)\right)^*=T_h(a^*)+o(1),
$$
where $a,b\in \Psi_{T^*X\times \mathbb{R}}(\Omega)$ are arbitrary elements and the estimate
$o(1)$ holds in the operator norm. This semiclassical quantization is a special case of
so-called \textit{asymptotic homomorphisms}, which play an important role in the theory of
$C^*$-algebras \cite{CoHi4}, \cite{Hig2}, \cite{Man1}. In particular, it follows that the
quantization $T_h$ induces the $K$-group homomorphism
$$
T:K_*(\overline{\Psi_{T^*X\times \mathbb{R}}(\Omega)}) \to K_*(\overline{\Sigma|_{\partial
M}}).
$$
Here we have used the fact that the algebras of operators with smooth symbols in question are
subalgebras in their closures and are closed with respect to holomorphic functional calculus.

Consider the commutative diagram
\begin{equation}\label{dia4}
\begin{array}{rcccccl}
0\to& \Psi^{-1}_{T^*X\times \mathbb{R}}(\Omega) & \to & \Psi_{T^*X\times \mathbb{R}}(\Omega) &
\to
& C^\infty(\partial S^*\mathcal{M}) & \to 0\\
& \;\;\;\downarrow T_h & & \downarrow T_h & & \downarrow t_h\\
 0\to &\Sigma_0 & \to & \Sigma|_{\partial M} & \to & C(\partial S^*\mathcal{M}) &
\to 0,
\end{array}
\end{equation}
where on the ideal $\Psi^{-1}_{T^*X\times \mathbb{R}}(\Omega)$ consisting of families of order
$\le -1$ we consider the restriction of the mapping $T_h$ and by $t_h$ we denote the
homomorphism on principal symbols.

The algebra $C^\infty(\partial M)$ is embedded in each of the algebras $\Psi_{T^*X\times
\mathbb{R}}(\Omega)$, $C^\infty(\partial S^*\mathcal{M})$, and $\Sigma|_{\partial M}$. The
diagram formed by the cones of these embeddings gives the square
\begin{equation}\label{quadrat}
\begin{array}{ccc}
K^*_c(\partial T^*\mathcal{M}) & \stackrel{\partial''}\longrightarrow & K^{*}_c(T^*X\times \mathbb{R})\\
\parallel & &\;\;\;\downarrow T \\
K^*_c(\partial T^*\mathcal{M}) & \stackrel{\partial'}\longrightarrow& K_{*}(\Sigma_0)
\end{array}
\end{equation}
of $K$-groups. (The mappings in this diagram are ordinary homomorphisms.) Here we have used the
isomorphism $K_*(\Psi^{-1}_{T^*X\times \mathbb{R}}(\Omega))\simeq K_*(C_0(T^*X\times
\mathbb{R}))$. The commutativity in \eqref{quadrat} follows from the fact that the boundary
mapping is natural with respect to asymptotic homomorphisms. (We leave the easy verification of
this fact to the reader; cf.~\cite{Hig2}.)

We embed the square \eqref{quadrat} in the diagram
$$
\begin{array}{ccc}
K_{*+1}(\partial M)  & \stackrel{\pi_*}\longrightarrow &  K_{*+1}(X)\\
\uparrow & & \uparrow \\
K^*_c(\partial T^*\mathcal{M}) & \stackrel{\partial''}\longrightarrow & K^{*}_c(T^*X\times \mathbb{R})\\
\parallel & &\;\;\;\downarrow T \\
K^*_c(\partial T^*\mathcal{M}) & \stackrel{\partial'}\longrightarrow& K_{*}(\Sigma_0)
\end{array}
$$
involving $K$-homology groups. The upper square in the diagram commutes. Indeed, the boundary
mapping $\partial''$ is given by the index of families with parameters ranging in $T^*X\times
\mathbb{R}$. This index can be computed by the index theorem in \cite{Nis1}:
$\partial''=\pi_!$, where
$$
\pi_!: K^*_c(T^*\partial M) \longrightarrow K_c^*(T^*X)
$$
is the direct image mapping in topological $K$-theory. It is also known \cite{BaDo1} that the
direct image $\pi_!$ gives the induced mapping $\pi_*$ after the passage to $K$-homology.

The commutativity of the square \fbox{A} in \eqref{dia3} will be proved once we establish that
the homomorphism $T$ is a $K$-group isomorphism inverse to the isomorphism $L$ described in the
preceding section.
\begin{lemma}
The mapping
$$
T: K_*(C_0(T^*X\times \mathbb{R}))\longrightarrow K_*(\Sigma_0)
$$
is the isomorphism and the inverse of $L$, see~\eqref{luke}.
\end{lemma}
\begin{proof}
To be definite, consider the mapping
$$
T: K_1(C_0(T^*X\times
\mathbb{R},\mathcal{K}))\longrightarrow K_1(\Sigma_0).
$$

Let us prove that the composition $LT$ gives the identity mapping of the $K$-group.

\textbf{1.} Indeed, on a symbol $u(\xi,p)$ invertible for $(\xi,p)\in T^*X\times \mathbb{R}$
and identically equal to unity outside a compact set, this mapping is defined by the formula
$$
 LT [u]=\ind T_hu\in
 K_c^0(T^*X),\quad [u]\in K_c^1(T^*X\times \mathbb{R}),
$$
for sufficiently small $h$. The index element is well defined, since the operator-valued
function$(T_h u)(\xi)$ has a compact fiber variation on $T^*X\setminus \mathbf{0}$ and is
invertible for $\xi\ne 0$.

\textbf{2.} By $\overline{T_h u}$ we denote the family of Fredholm operators on $T^*X\setminus
\mathbf{0}$ coinciding with the family $T_h u$ for $|\xi|<1$ and equal to the family
$$
u\BL((\stackrel 2 r+|\xi|-1)\xi,ih \stackrel  1 {r\frac\partial{\partial r}}+ih\frac{n+1}2\BR)
$$
for $|\xi|\ge 1$. For sufficiently small $h$, this family will be invertible for all $\xi$
(this follows from the boundedness of the support of $1-u$ and also from the fact that the
estimates \eqref{7} hold uniformly with respect to the parameter $\la$, $\const>\lambda\ge 0$,
for the symbol $u(\xi(e^{-t}+\lambda),p)$). By construction, we have
$$
\ind T_h u=\ind  \overline{T_h u}.
$$
But the family $\overline{T_h u}$ consists of identity operators for $\xi$ outside a compact
set. Hence its index can be computed by formula \eqref{indf1} and is equal to the index of the
family of conormal symbols $u(\xi,p)$ (modulo periodicity), i.e., indeed gives the original
element
$$
LT [u]=[u]
$$
of the $K$-group.
\end{proof}

\begin{remark}
The introduction of the semiclassical parameter $h\to 0$ can also be viewed as an adiabatic
limit (e.g., see~\cite{Che3}) reducing the study of edge symbols as operators on the infinite
cone $K_\Omega$ to families with parameters of operators on sections of the cone.  Note that
families with parameters play a role similar to that of the edge symbol in the theory of
elliptic operators on manifolds with fibered cusps~\cite{MaMe3}. It would be of interest to
clarify the relationship not only between symbols, but also between operators in these two
theories.
\end{remark}

\section{Some Remarks}

\paragraph{1. The classification of edge morphisms.}
Let us show that the classification of elliptic edge morphisms with conditions and
co-conditions on the edge (e.g., see~\cite{NSScS98}) follows from the classification of edge
operators. By $\widetilde{\Sigma}\supset \Sigma$ we denote the symbol algebra corresponding to
the algebra of zero-order pseudodifferential operators with edge and co-edge conditions and
with pseudodifferential operators on the edge.

\begin{lemma}\label{lem51}
The embedding $\Con(C(M)\to \Sigma)\subset \Con(C(M)\to \widetilde{\Sigma})$ induces an
isomorphism in $K$-theory. The embedding $\Sigma_0\subset \widetilde{\Sigma}_0$ has a similar
property.
\end{lemma}
\begin{proof}
The first fact follows from the commutative diagram
$$
\begin{array}{ccccc}
SC(S^*X,\mathcal{K}) & \to& \Con(C(M)\to \Sigma) &
\stackrel{(\sigma,\sigma_c)}\longrightarrow &    \im (\sigma,\sigma_c)\\
\cap & & \cap & & \parallel \\
SC(S^*X,\widetilde{\mathcal{K}}) & \to& \Con(C(M)\to \widetilde{\Sigma}) &
\stackrel{(\sigma,\sigma_c)}\longrightarrow & \im (\sigma,\sigma_c)
\end{array}
$$
(the projections correspond to taking the interior and conormal symbols), since the left
vertical embedding induces an isomorphism of $K$-groups. Here by $\widetilde{\mathcal{K}}$ we
denote the algebra of compact operators in the direct sum $\mathcal{K}^{0,0}(K_\Omega)\oplus
\mathbb{C}$, and the embedding is induced by the representation of operators in $\mathcal{K}$
as the upper left corner of the matrix.

The second isomorphism $K_*(\Sigma_0)\simeq K_*(\widetilde{\Sigma}_0)$ can be obtained in a
similar way.
\end{proof}

We set  $\widetilde{A}=\Con(C(M)\oplus C(X)\to \widetilde{\Sigma})$. The group
$K_0(\widetilde{A})$ classifies the stable homotopy classes of elliptic \emph{edge problems}
(see~\cite{NSScS98})
$$
\mathcal{D}=\left(\begin{array}{cc} D & C \\ B & P
\end{array}\right):\mathcal{W}^{0,0}(\mathcal{M},E)\oplus L^2(X,E_0)\longrightarrow
\mathcal{W}^{0,0}(\mathcal{M},F)\oplus L^2(X,F_0),
$$
$E,F\in \Vect(M)$, $E_0,F_0\in \Vect(X)$.

\begin{proposition}
One has the isomorphism $K_0(\widetilde{A})\simeq K_0(A)\oplus K^0(X)$.
\end{proposition}
\begin{proof}
Consider the sequence
$$
0\to \Con(C(M)\to\widetilde{\Sigma})\subset \Con(C(M)\oplus C(X)\to \widetilde{\Sigma})\to
C(X)\to 0.
$$
The corresponding sequence
\begin{equation}
\label{kgroups} \ldots\to K_*(A) \longrightarrow K_*(\widetilde{A}) \longrightarrow
K^*(X)\to\ldots
\end{equation}
of $K$-groups (here we have used Lemma~\ref{lem51}) splits.

First, let us indicate a splitting $j:K^0(X)\to K_0(\widetilde{A})$. To this end, we choose an
edge symbol
$$
\sigma_\Lambda^0(x):\cK^{0,0}(K_{\Omega_x})\longrightarrow \cK^{0,0}(K_{\Omega_x})
$$
such that for all $(x,\xi)\in T^*X\setminus\mathbf{0}$ it is Fredholm, has a one-dimensional
kernel, a trivial cokernel, and unit principal symbol (e.g., see~\cite{NSScS3}).

Now we define an element $j[E]$ for each vector bundle $E\in\Vect(X)$ as the following symbol:
$$
a=1, a_\Lambda=\left(
                 \begin{array}{c}
                    \sigma_\Lambda^0\otimes 1_E \oplus 1\\
                    i^*
                 \end{array}
               \right):\cK^{0,0}(K_\Omega,E\oplus E^\perp)\longrightarrow
               \begin{array}{c}
                  \cK^{0,0}(K_\Omega,E\oplus E^\perp)\\
                  \oplus \\
                  E
               \end{array},
$$
where by $i:E\to \mathcal{K}^{0,0}(K_\Omega)$ we denote some embedding of the
finite-dimensional bundle into the infinite-dimensional bundle, $i^*$ is the adjoint mapping,
and $E^\perp$ is the complementary bundle.

A similar splitting in the term $K_1(\widetilde{A})$ can be obtained with the use of suspension
($K^1(X)\simeq K(X\times \mathbb{S}^1)/K(X)$).
\end{proof}

\paragraph{2. Topological obstructions in edge theory.}
The commutativity of the square \fbox{A} in diagram~\eqref{dia3} can be interpreted as a
topological obstruction.
\begin{corollary}[\cite{NSScS4}]
The stable homotopy class of an elliptic interior symbol $\sigma$ on $T^*\mathcal{M}$ contains
a representative possessing a compatible elliptic edge symbol (i.e., there exists a Fredholm
operator) if and only if
$$
\pi_![\sigma|_{\partial T^*\mathcal{M}}]=0\in K^1_c(T^*X), \qquad \pi:\partial M\longrightarrow
X,
$$
where $\pi_!:K^0_c(\partial T^*\mathcal{M})\to K^1_c(T^*X)$ is the direct image mapping.
\end{corollary}
A similar obstruction for self-adjoint elliptic operators is given by the same formula. We
point out that there is an essential difference between the even and odd cases.  Namely, the
value of the boundary mapping $K^0_c(T^*\mathcal{M}) \stackrel \partial\to K_1(S\Sigma_0
)\simeq K^1_c(T^*X)$ on the element determined by an elliptic symbol $\sigma$ is expressed in
terms of the index
$$
\partial [\sigma]=p \ind \sigma_\Lambda\in K^0(S^*X)/K^0(X)\simeq K^1_c(T^*X),
$$
where $p:K(S^*X)\to K(S^*X)/ K(X)$, of a compatible edge symbol $\sigma$. (The edge symbol is
Fredholm, since by~\cite{NSScS5} for an elliptic interior symbol $\sigma$ there always exists a
compatible invertible conormal symbol.)

However, for self-adjoint elliptic operators there does not necessarily exist a compatible
elliptic self-adjoint conormal symbol. This is shown by the following example even in the case
of conical singularities.
\begin{example}
Consider an interior elliptic symmetric operator of Dirac type on a manifold with $N$ conical
points of the form
$$
D=\Gamma\left(\frac \partial {\partial t}+A\right), \qquad r=e^{-t}
$$
in a neighborhood of the singularities, where $\Gamma^2=-1$, $\Gamma A+A\Gamma=0$, and
$\Gamma^*=-\Gamma, A=A^*$. Then the boundary mapping $K^1_c(T^*\mathcal{M})\to
K^0_c(T^*X)=\mathbb{Z}^N$, whose triviality on the element $[\sigma(D)]$ is the obstruction to
the existence of a self-adjoint elliptic operator, takes $[\sigma(D)]$ to the sequence of
indices of operators induced on the bases of the cones:
$$
\ind \bigl( A:\im (\Gamma+i)\longrightarrow  \im(\Gamma-i)\bigr).
$$
Note that the sum of these indices is zero, since the index is cobordism invariant \cite{Pal1},
but separate terms can be nonzero. For example, one can take $M=[0,1]\times \mathbb{CP}^{2n}$
with the signature operator.
\end{example}

\section{Appendix}

\paragraph{1. Semiclassical quantization.}
Let $\Om$ be a smooth closed manifold. By
$\Psi^0\equiv\Psi^0(\Om,\RR^{l+1}_{\eta,p})\phantom{\Bigr.}$ we denote the algebra of classical
pseudodifferential operator of order zero with parameters $(\eta,p)\in\RR^{l+1}$ in the sense
of Agranovich--Vishik on $\Om$, and by $\Psiw\equiv\Psiw(K_\Om)$ we denote the algebra of edge
symbols\footnote{Edge symbols depend also on $x$, but we omit this dependence, since it does
not affect the construction in any way.} of order zero depending on the parameter $\xi\in
\mathbb{S}^{l-1}\subset\RR^l_\xi$ in the space $\mathcal{K}^{0,0}(K_\Omega)$ on the cone
$K_\Om$.

The problem is to construct an asymptotic homomorphism (quantization)
\begin{equation*}
    T_h:\Psi^0\lra\Psiw,\quad h\in (0,1],
\end{equation*}
such that the interior symbol and the conormal symbol of the edge symbol given by the
quantization satisfy the relations
\begin{align}
    \si(T_h a)(\eta,p)&=\si(a)(q,\eta,hp),\label{1}\\
    \si_c(T_h a)(\eta,p)&=a\left(0,h\left(p+i(n+1)/2\right)\right).\label{2}
\end{align}

\begin{theorem}
The mapping $T_h$ can be defined by the formula
\begin{equation}\label{3}
    (T_h a)(\xi)=
    a\BL(\xi \ovs2r,ih\ovs1{r\pd{}r}+ih\frac{n+1}2\BR),
\end{equation}
where numbers over operators are treated in the sense of noncommutative analysis
\rom{\cite{NSS8}} and functions of $ir\pa/\pa r$ are defined with the help of the Mellin
transform on the weight line $\im p=-(n+1)/2$, where $n=\dim\Omega$.
\end{theorem}

\begin{proof}
First, let us show that the mapping~\eqref{3} is well defined and continuous in the spaces
\begin{equation*}
    T_h:\Psi^0\lra \cB(\cK^{0,0}(K_\Om))
\end{equation*}
and indeed specifies an asymptotic quantization, i.e., satisfies
\begin{equation*}
    T_h(a)T_h(b)=T_h(ab)+o(1),\quad h\to0,
\end{equation*}
where the estimate $o(1)$ on the right-hand side is in the operator norm. To prove this, we make
the change of variables $r=e^{-t}$. It takes the space $\cK^{0,0}(K_\Om)$ to the space
$L^2(\RR_t,L^2(\Om),e^{(n+1)t}dt)$ of functions ranging in $L^2(\Om)$ and square integrable with
respect to the measure $e^{(n+1)t}dt$ on the real line, and the quantization formula~\eqref{3}
acquires the standard $h$-pseudodifferential form
\begin{equation}\label{4}
    T_h(a(\eta,p))=
    e^{(n+1)t/2}A\BL(\ovs2t,-ih\ovs1{\pd{}t}\BR)e^{-(n+1)t/2},
\end{equation}
where
\begin{equation}\label{5}
    A(t,p)=a(\xi e^{-t},p).
\end{equation}
Note that hence it suffices to study the operator
$$
A\BL(\ovs2t,-ih\ovs1{\pd{}t}\BR):L^2(\RR_t\times\Om,dtd\om)\lra L^2(\RR_t\times\Om,dtd\om)
$$
in the standard Lebesgue space $L^2(\RR\times\Om)=L^2(\RR_t\times\Om,dtd\om)$. Since the
operator-valued function $a(\eta,p)$ is a zero-order pseudodifferential operator with parameters
in the sense of Agranovich--Vishik, it follows that its norm in $L^2(\Om)$ satisfies the
estimates
\begin{equation}\label{6}
    \norm{\pd{^{\a+k}a(\eta,p)}{\eta^\a\pa p^k}}\le C_{\a k}(1+\abs{\eta}+\abs{p})^{-\abs{\a}-k}
\end{equation}
for $\abs{\a}+k=0,1,2,\ldots$ and hence the norm of the function~\eqref{5} satisfies the
estimates
\begin{equation}\label{7}
    \norm{\pd{^{l+k}A(t,p)}{t^l\pa p^k}}\le \wt C_{lk}(1+\abs{p})^{-k}.
\end{equation}
Indeed,
\begin{equation*}
    \pd{^{l+k}A(t,p)}{t^l\pa p^k}=\sum_{\abs{\a}\le l}M_{l\a}e^{-\abs{\a}t}
    \pd{^{\a+k}a}{\eta^\a\pa p^k}(\xi e^{-t},p),
\end{equation*}
where the coefficients $M_{l\a}$ depend only on $\xi$, and hence (we assume that $\xi$ lies on
the sphere)
\begin{equation*}
    \norm{\pd{^{l+k}A(t,p)}{t^l\pa p^k}}\le \sum_{\abs{\a}\le l}\abs{M_{l\a}}C_{\a k}
    e^{-\abs{\a}t}(1+e^{-t}+\abs{p})^{-\abs{\a}-k}.
\end{equation*}
It remains to note that
\begin{equation*}
    e^{-t}(1+e^{-t}+\abs{p})^{-1}\le1,\quad (1+e^{-t}+\abs{p})^{-1}\le(1+\abs{p})^{-1}.
\end{equation*}

Hence we have shown that the operator-value symbol $A(t,p)$ belongs to the H\"ormander class
$S^0\equiv S^0_{1,0}$.

\begin{lemma}\label{nka}
The following assertions hold.

\rom{(i)} If $A\in S^0$, then the operator
\begin{equation*}
    \wh A\equiv A\BL(\ovs2t,-ih\ovs1{\pd{}t}\BR):L^2(\RR\times\Om)\lra L^2(\RR\times\Om)
\end{equation*}
is bounded uniformly with respect to $h\in(0,1]$.

\rom{(ii)} If $A,B\in S^0$, then
\begin{equation*}
    \wh{AB}-\wh A\wh B=O(h).
\end{equation*}
\end{lemma}

\begin{proof}
The proof of this lemma is standard.  It follows the proof of Theorem IV.6
in~\cite[p.~282]{NSS8} almost word for word, and we omit it.
\end{proof}

Thus we have shown that $T$ is an asymptotic homomorphism. It remains to verify that the
interior and conormal symbols of the operator $T(a)$ satisfy relations~\eqref{1} and \eqref{2}.
But this follows directly from Definition~\ref{vot-vam}.
\end{proof}

\paragraph{2. Estimates modulo compact operators.}

\begin{theorem}
\rom{(i)} The norm modulo compact operators of an edge-degenerate zero-order pseudodifferential
operator in the space$\mathcal{W}^{0,0}(\cM)$ is equal to
\begin{equation}
\inf_{K\in \mathcal{K}}\|D+K\|_{\mathcal{B}(\cW^{0,0}(\mathcal{M})
)}=\max\left(\max_{S^*\mathcal{M}} |\sigma(D)|,\max_{S^*X}
\|\sigma_\Lambda(D)\|_{\mathcal{B}(\cK^{0,0}(K_\Omega))}\right)
\end{equation}

\rom{(ii)} The norm modulo compact operators of an edge symbol is equal to
\begin{multline}
\inf_{K\in C(S^*X,\mathcal{K})}\max_{S^*X}\|a_\Lambda+K\|_{\mathcal{B}
(\cK^{0,0}(K_\Omega))}\\= \max\left(\max_{\partial S^*\mathcal{M}}
|\sigma(a_\Lambda)|,\sup_{X\times \mathbb{R}}
\|\sigma_c(a_\Lambda)\|_{\mathcal{B}(L^2(\Omega))}\right)
\end{multline}
\end{theorem}

\begin{proof}
The proof of both assertions follows the same scheme and is very simple. The symbol (interior
or edge) of the operator $D$ can be computed (up to Fourier transform with respect to part of
variables, which does not change the norm) as the strong limit of the expression
$U_\la^{-1}AU_\la$, where $U_\la$ is a unitary local dilation group in a neighborhood of the
corresponding point of $\cM$. It follows that the norm of the symbol does not exceed the norm
of the operator. Next, the product of operators corresponds to the product of symbols, and the
adjoint operator corresponds to the adjoint symbol. It follows that the mapping ``operator
$\mapsto$ symbols'' extends to a $*$-homomorphism of the $C^*$-algebra obtained by the closure
of the algebra of pseudodifferential operator into the corresponding $C^*$-algebra of symbols.
The kernel of this homomorphism coincides with the set of compact operators, since the equality
of all symbols to zero is a necessary and sufficient condition for compactness. By the general
properties of $*$-homomorphisms, we see that the mapping ``operator $\mapsto$ symbols'' is an
isometric isomorphism of the quotient algebra of pseudodifferential operators modulo compact
operators onto the symbol algebra, which proves the desired assertion.

\smallskip

To prove the corresponding assertion for edge symbols, we should write out a local
one-parameter group computing the interior symbol of the edge symbol. This (additive) group
acts in a neighborhood of the point $(r,\om_0)$ of the cone $K_\Om$ in the coordinates
$(r,\tau)=(r,r(\om-\om_0))$ by the formula
\begin{equation}\label{group}
    U_\la\ph(r,\tau)=\ph(r+\la,\tau).
\end{equation}
The conormal symbol $a(0,p)$ (more precisely, its Mellin transform) can be computed as the
strong limit
\begin{equation*}
    \si_c(T(a))=s\text{-}\lim_{\la\to\infty} U_\la^{-1}T(a)U_\la,
\end{equation*}
where $U_\la$ is the group of dilations with respect to the variable $r$.
\end{proof}

%\bibliographystyle{unsrt1}
%\bibliography{elliptic}
%\end{document}

\end{document}